%% file: m3-I-8.tex
\input m3-macs

\pageno=81

\tinfo I.8.81-89

\SetTFLinebox{\gtp }
\SetFLinebox{\gtv3 }
\SetHLinebox{\issn}

\H 8. Explicit formulas for the Hilbert symbol 

Sergei V. Vostokov

\SetAuthorHead{S. Vostokov}
\SetTitleHead{Part I. Section 8. Explicit formulas for the Hilbert symbol \qquad\qquad}

Recall that
the Hilbert symbol for
a local field $K$ with finite residue field
which contains a primitive $p^n$th
root of unity $\zeta_{p^n}$ is a pairing
$$
(\,\,,\,\,)_{p^n}\colon K^*/K^{*\, p^n}\times K^*/K^{*\, p^n}\to \langle \zeta_{p^n}\rangle,
\quad (\alpha,\beta)_{p^n}=
\gamma^{\Psi _{K}(\alpha )-1}, \quad \gamma^{p^n}=\beta, $$ 
where 
$
\Psi_{K}\colon  K^* \to \Gal(K^{\ab}/K)
$
is the reciprocity map.

\HH 8.1. History of explicit formulas for the Hilbert symbol

 There are two different branches of explicit reciprocity formulas (for the
Hilbert symbol).

\HHH 8.1.1. The first branch (Kummer's type formulas)

\th Theorem {{\rm (E. Kummer 1858)}} 

Let $K=\Bbb Q_p(\zeta_p)$, $p\not=2$.
Then for principal units $\varepsilon ,\eta$
$$
\boxed{(\varepsilon,\eta)_p=\zeta_p^{\dsize \,\res (\lln \eta(X)\,d\lln \varepsilon(X) X^{-p})} }
$$
where $\varepsilon(X)|_{X=\zeta_p-1}=\varepsilon$, 
$\eta(X)|_{X=\zeta_p-1}=\eta$,
$\varepsilon(X), \eta(X)\in \Bbb Z_p[[X]]^*$.
\endth

The important point is that one associates to the elements 
$\varepsilon, \eta$ 
the series 
$\varepsilon(X), \eta(X)$ 
in order to calculate the value of the Hilbert symbol.

\th Theorem {{\rm (I. Shafarevich 1950)}}

Complete explicit formula for the Hilbert norm residue symbol
$(\alpha,\beta)_{p^n}$, $\alpha,\beta\in K^*$, 
$K\supset \Bbb Q_p(\zeta_{p^n})$, $p\not=2$, 
using a special basis of the group of principal units.
\endth

This formula is not very easy to use
because of the special  basis of the group of units
and certain difficulties with its verification for $n>1$. 
One of applications of this formula was in the work of Yakovlev
on the description of the absolute Galois group of a local
field in terms of generators and relations.

Complete formulas, which are simpler that Shafarevich's formula,
were discovered in the seventies:

\th Theorem {{\rm (S.~Vostokov 1978), (H. Br\"uckner 1979)}}

Let a local field $K$ with finite residue field contain $\Bbb Q_p(\zeta_{p^n})$ and let $p\not=2$.
Denote $\coo =W(k_K)$, $\Tr=\Tr_{\coo /\Bbb Z_p}$.
Then for $\alpha,\beta\in K^*$
$$
\boxed{ 
(\alpha,\beta)_{p^n}=\zeta_{p^n}^{\dsize \,\Tr\res \,\Phi(\alpha,\beta)/\Underline{s}}, \quad 
\Phi(\alpha,\beta)=
l(\Underline{\beta})\Underline{\alpha}^{-1}d\Underline{\alpha}-l(\Underline{\alpha})\frac{1}{p}\Underline{\beta}^{-\triangle}
d\Underline{\beta}^{\triangle} }  $$
where
$\Underline{\alpha}=\theta X^m(1+\psi(X))$, $\theta\in \Cal R$, $\psi\in X\coo [[X]]$, 
is such that $\Underline{\alpha}(\pi)=\alpha$,
$\Underline{s}=\Underline{\zeta_{p^n}}^{p^n}-1$,
$$
\aligned  
l(\Underline{\alpha})&=\frac{1}{p}\log\bigl(\Underline{\alpha}^p/\Underline{\alpha}^{\triangle}\bigr),\\
\left (\sum a_iX^i\right)^{\triangle} &=\sum \Frob_K(a_i)X^{pi} ,\quad a_i\in\coo .
\endaligned
$$
\endth

Note that for the term $X^{-p}$ in Kummer's theorem
can be written as 
$X^{-p}=1/(\Underline{\zeta_p}^p-1) \mod p$, 
since $\zeta_p =1+ \pi$ and so 
$\Underline{s}=\Underline{\zeta_p}^p-1=(1+X)^p-1=X^p \mod p$.

The works \cite{V1} and \cite{V2}  contain two different proofs of this formula.
One of them is to construct  the explicit pairing
$$(\alpha,\beta)\to \zeta_{p^n}^{\dsize \,\Tr\res \,\Phi(\alpha,\beta)/\Underline{s}}
$$ and check the correctness of the definition and all the properties
of this pairing  completely independently of class field theory
(somewhat similarly to how one works with the tame symbol), 
and only at the last step to show that the pairing 
coincides with the Hilbert symbol.
The second method, also followed by  Br\"ukner, is different:
it uses Kneser's (1951) calculation 
of symbols and reduces the problem to a simpler one:
to find a formula for 
$(\varepsilon,\pi)_{p^n}$ where $\pi$ is a prime element of $K$
and $\varepsilon$ is a principal unit of $K$.
Whereas the first method is very universal and can be extended to formal groups and higher local fields, the second method works well in the classical situation
only. 

For $p=2$ explicit formulas were obtained by
G. Henniart (1981) who followed to a certain extent  Br\"uckner's method,
and S.~Vostokov and I. Fesenko (1982, 1985).

\HHH 8.1.2. The second branch (Artin--Hasse's type formulas)

\th Theorem {{\rm (E. Artin and H. Hasse 1928)}}

Let $K=\Bbb Q_p(\zeta_{p^n})$, $p\not=2$.
Then for a principal unit $\varepsilon $
and prime element $\pi=\zeta_{p^n}-1$ of $K$
$$ 
\boxed{
(\varepsilon,\zeta_{p^n})_{p^n}=\zeta_{p^n}^{\dsize\, \Tr(-\lln\varepsilon)/p^n}, \quad 
(\varepsilon,\pi)_{p^n}=\zeta_{p^n}^{\dsize \, \Tr(\pi^{-1}\zeta_{p^n}\lln\varepsilon)/p^n}
}
$$
where $\Tr=\Tr_{K/\Bbb Q_p}$.
\endth

\th Theorem {{\rm (K. Iwasawa 1968)}}

Formula for $(\varepsilon,\eta)_{p^n}$
where $K=\Bbb Q_p(\zeta_{p^n})$, $p\not=2$,
$\varepsilon,\eta$ are principal units of $K$
and  $v_K(\eta-1)>2v_K(p)/(p-1)$.
\endth

To some extent the following formula can be viewed
as a formula of Artin--Hasse's type.
 Sen deduced it using his theory of continuous Galois representations
which itself is a  generalization
of a part of Tate's theory of $p$-divisible groups.
The Hilbert symbol is  interpreted as the cup product of $H^1$.

\th  Theorem {{\rm (Sh. Sen 1980)}}

Let $| K:\Bbb Q_{p}| <\infty $, $\zeta _{p^{n}}\in K$,
and let  
$\pi $ be a prime element of $\Cal O_{K}$.
Let $g(T),h(T)\in W(k_{K})[T]$
be such that $g(\pi )=\beta \not=0$, $h(\pi )=\zeta _{p^{m}}$.  
Let $\alpha \in 
\Cal O_{K}$, $v_{K}(\alpha )\ge 2v_K(p)/(p-1)$.
Then  
$$
\boxed{
(\alpha ,\beta )_{p^{m}}=\zeta _{p^{m}}^{c},\quad c=\frac{1}{%
p^{m}}\Tr_{K/\Bbb Q_{p}}\biggl( \frac{\zeta_{p^{m}}}{h^{\prime }(\pi )%
}\frac{g^{\prime }(\pi )}{\beta }\lln \alpha \biggr).
}
$$
\endth

\smallskip

R. Coleman  (1981) gave a new form of explicit formulas
which he proved for  $K=\Bbb Q_p(\zeta_{p^n})$. 
He uses formal power series associated to norm compatible sequences
of elements in the tower of finite subextensions
of the $p$-cyclotomic extension of the ground field 
and his formula can be viewed as a generalization of Iwasawa's formula.

\HH 8.2. History: 
Further developments

\HHH 8.2.1

Explicit formulas for the (generalized) Hilbert symbol
in the case where it is defined by an appropriate class field theory.

\df Definition

Let $K$ be an $n$-dimensional local field
of characteristic 0 which contains a primitive $p^m$th
root of unity.
The  $p^{m}$th Hilbert
symbol is defined as 
$$
K_{n}^{\tpp}(K)/p^{m}\times K^{* }/K^{* p^{m}}\rightarrow \langle
\zeta_{p^m} \rangle, \quad 
(\alpha ,\beta )_{p^m}=\gamma ^{\Psi _{K}(\alpha )-1},
\quad \gamma^{p^m}=\beta, 
$$
where 
$
\Psi _{K}\colon K_{n}^{\tpp}(K)\to \Gal(K^{\ab}/K)
$
is the reciprocity map. 
\enddf

For higher local fields and $p>2$ complete formulas of Kummer's type 
were constructed by S.~Vostokov (1985).
They are discussed in subsections 8.3 and 
their applications to $K$-theory of higher local fields
and $p$-part of the existence theorem 
in characteristic 0 are discussed in
subsections 6.6, 6.7 and 10.5. 
For higher local fields, $p>2$ and  Lubin--Tate formal group
complete formulas of Kummer's type 
were deduced by I. Fesenko (1987).

Relations of the formulas with syntomic cohomologies
were studied by K.~Kato (1991) in a very important work 
where it is suggested to use Fontaine--Messing's syntomic cohomologies 
and an interpretation of the Hilbert symbol as the cup product
explicitly computable in terms of the cup product of syntomic
cohomologies; this approach implies Vostokov's formula.
On the other hand, Vostokov's formula appropriately generalized
defines a homomorphism from the Milnor $K$-groups to
cohomology groups of a syntomic complex (see subsection~15.1.1). 
M. Kurihara (1990) applied syntomic cohomologies to  deduce 
Iwasawa's and Coleman's formulas in the multiplicative case. 

For higher local fields complete formulas of Artin--Hasse's type were 
constructed by
M. Kurihara (1998), see section 9.

\HHH 8.2.2. Explicit formulas for $p$-divisible groups 

\df Definition

Let $F$ be a formal $p$-divisible group over the ring $\Cal O_{K_0}$
where $K_0$ is a subfield of a local field $K$.
Let $K$ contain $p^n$-division points of $F$.
Define the Hilbert symbol 
by
$$ K^*\times F(\Cal M_K)\to \kr [p^n],\quad 
(\alpha,\beta)_{p^n}=
\Psi_{K}(\alpha)(\gamma)-_F\gamma, \quad [p^n](\gamma)=\beta, $$ 
where 
$
\Psi_{K}\colon  K^* \to \Gal(K^{\ab}/K)
$
is the reciprocity map. 
\enddf

For  formal Lubin--Tate groups,
complete formulas of Kummer's type
were obtained by S.~Vostokov (1979) for odd $p$
and S. Vostokov and I. Fesenko (1983) for even $p$.
For  relative formal Lubin--Tate groups
complete formulas of Kummer's type
were obtained by S.~Vostokov and A. Demchenko (1995).

For  local fields with finite residue field and formal Lubin--Tate groups
formulas of Artin--Hasse's type
were deduced by  A. Wiles (1978) for
$K$ equal to  the $[\pi^n]$-division field of the isogeny $[\pi]$
of a formal Lubin--Tate group;
by V. Kolyvagin (1979)
for $K$ containing the $[\pi^n]$-division field of the isogeny $[\pi]$;
by R. Coleman (1981) in the multiplicative case and 
some partial cases of Lubin--Tate groups;
his conjectural formula in the general case of Lubin--Tate groups was proved by
 E. de Shalit (1986) 
for $K$ containing the $[\pi^n]$-division field of the isogeny $[\pi]$.
This formula was generalized by Y. Sueyoshi (1990) for relative formal Lubin--Tate groups.
F. Destrempes (1995) extended Sen's formulas to Lubin--Tate formal groups. 

J.--M. Fontaine (1991) used his crystalline ring
and his and J.--P. Wintenberger's theory of field of norms
for the $p$-cyclotomic extension to relate 
Kummer theory with Artin--Schreier--Witt theory
and deduce in particular some formulas of Iwasawa's type
using Coleman's power series. 
D. Benois (1998) further extended this approach 
by using Fontaine--Herr's complex and deduced Coleman's formula. 
V. Abrashkin (1997) used another arithmetically profinite extension
($L=\cup F_i$ of $F$, $F_i=F_{i-1}(\pi_i)$, $\pi_i^p=\pi_{i-1}$, $\pi_0$ being
a prime element of $F$)
to deduce the formula of Br\"uckner--Vostokov.  

For formal groups
which are defined over an  absolutely unramified local field $K_0$ 
($e(K_0|\Bbb Q_p)=1$) 
 and therefore are parametrized by Honda's systems,
formulas of Kummer's type
were deduced by D. Benois and S.~Vostokov (1990), for $n=1$
and one-dimensional formal groups, and 
by V. Abrashkin (1997) for arbitrary $n$ and arbitrary formal group 
with  restriction that $K$ contains a primitive $p^n$th root of unity.
For one dimensional formal groups and arbitrary $n$ 
without restriction that $K$ contains a primitive $p^n$th root of unity 
in the ramified case
formulas were obtained by S.~Vostokov and A. Demchenko (2000).
For arbitrary $n$ and arbitrary formal group without restrictions on $K$
Abrashkin's formula was established by Benois (2000),  
see subsection 6.6 of Part~II.

Sen's formulas were generalized to all  $p$-divisible groups 
by D. Benois (1997) using 
an interpretation of the Hilbert pairing in terms of 
an explicit construction of $p$-adic periods. 
T. Fukaya (1998) generalized the latter for higher local fields.

\HHH 8.2.3. Explicit formulas for $p$-adic representations 

The previously discussed explicit formulas can be viewed 
as  a description of the
exponential map from the tangent space of a formal group 
to the first cohomology group with coefficients in the Tate module. 
Bloch and Kato (1990) defined a generalization of the 
exponential map to de Rham representations. An explicit 
description of this map is closely related to the
computation of Tamagawa numbers of motives which play an important role
in  the Bloch--Kato conjecture.  The description of this map for
the $\Bbb Q_p\,(n)$ over cyclotomic fields was given by Bloch--Kato (1990)
and Kato (1993); it can be viewed as a vast generalization of Iwasawa's
formula (the case $n=1$).
B. Perrin-Riou constructed
an Iwasawa theory for  crystalline representations 
over an absolutely unramified local field 
and conjectured an explicit description of the cup product 
of the cohomology groups. 
There are three different approaches which culminate in the proof
of this conjecture by P. Colmez (1998),
K.~Kato--M. Kurihara--T. Tsuji (unpublished) and 
for crystalline representations of finite height
by D. Benois (1998).

K.~Kato (1999)
gave generalizations of explicit formulas of Artin--Hasse, Iwasawa and Wiles type to $p$-adically complete discrete valuation fields 
and $p$-divisible groups which relates norm compatible sequences
in the Milnor $K$-groups and trace compatible sequences in
differential forms;
these formulas are applied in his other work
to give an explicit description 
in  the case of $p$-adic completions of function fields of modular curves.

\HH 8.3. Explicit formulas in higher dimensional fields of characteristic 0

Let $K$ be an $n$-dimensional 
field of characteristic 0, $\chr(K_{n-1})=p$, $p>2$. 
Let  $\zeta _{p^{m}}\in K$.

Let $t_1,\dots,t_n$ be a system of local parameters of $K$.

For an element $$\alpha=t_n^{i_n}\dots t_1^{i_1}\theta
(1+\sum a_{J} t_n^{j_n}\dots t_1^{j_1}), \quad \theta\in\Cal R^*, a_J\in W(K_0), $$
$(j_1,\dots,j_n)>(0,\dots,0)$ 
denote by 
$\Underline{\alpha}$ the following element
$$X_n^{i_n}\dots X_1^{i_1}\theta
(1+\sum  a_{J} X_n^{j_n}\dots X_1^{j_1})$$
in $F\{\!\{X_1\}\!\}\dots\{\!\{X_n\}\!\}$
where $F$ is the fraction field of $W(K_0)$. 
Clearly $\Underline{\alpha}$ is not uniquely determined
even if the choice of a system of local parameters is fixed.

{\it Independently} of class field theory define
the following explicit map
$$V(\,\,,\,)_{m}\colon (K^{* })^{n+1}\to \langle\zeta _{p^{m}} \rangle $$
by the formula 
$$
\aligned 
&V(\alpha _{1},\dots ,\alpha _{n+1})_{m}=\zeta _{p^{m}}
^{\dsize \, \Tr\res\,\,\Phi
(\alpha _{1},\dots ,\alpha _{n+1})/\Underline{s}}, \quad
\Phi (\alpha _{1},\dots ,\alpha _{n+1}) \\ 
&\quad =
\sum_{i=1}^{n+1}\frac{(-1)^{n-i+1}%
}{p^{n-i+1}}l\left( \Underline{\alpha _{i}}\right) \frac{d\Underline{\alpha _{1}}}{\Underline{\alpha
_{1}}}\wedge \cdots \wedge \frac{d\Underline{\alpha _{i-1}}}{
\Underline{\alpha _{i-1}}} 
\wedge \frac{d\Underline{\alpha _{i+1}}^{\triangle }}{\Underline{\alpha _{i+1}}^{\triangle }}%
\wedge \cdots \wedge \frac{d\Underline{\alpha _{n+1}}^{\triangle }}{
\Underline{\alpha
_{n+1}}^{\triangle}}
\endaligned 
$$ 
where $\Underline{s}=\Underline{\zeta_{p^{m}}}^{p^{m}}-1$, $\Tr=\Tr_{W(K_{0})/\Bbb Z%
_{p}}$, $\res=\res_{X_1,\dots, X_n}$, 
$$ 
l\left( \Underline{\alpha}\right) =\frac{1}{p}\lln \left( \Underline{\alpha}^{p}/\Underline{\alpha}^{\triangle }\right), \quad 
\bigl( \sum a_{J}X_{n}^{j_{n}}\cdots X_{1}^{j_{i}}\bigr) ^{\triangle
}=\sum \Frob (a_{J})X_{n}^{pj_{n}}\cdots X_{1}^{pj_{1}}.
$$

\th Theorem 1

The map $V(\,\,,\,)_m$ is well defined,
 multilinear and symbolic. It induces a homomorphism
$$
K_{n}(K)/p^{m}\times K^{* }/K^{*\,p^{m}}\to \mu _{p^{m}}
$$
and since $V$ is sequentially continuous, a homomorphism 
$$
V(\,\,,\,)_{m}\colon K_{n}^{\tpp}(K)/p^{m}\times K^{* }/K^{*\,p^{m}}\to \mu
_{p^{m}}
$$
which is non-degenerate.
\endth

\pf Comment on Proof

A set of elements $t_{1},\dots ,t_{n}$, $%
\varepsilon _{{\bold j}},\omega $ (where $\bold j$ runs over
a subset of $\Bbb Z^n$) is called a {\it Shafarevich basis}  
of $K^*/K^{*p^m}$ if
\Roster
\Item{(1)}
every $\alpha \in K^{* }$ can be written
as a convergent product 
$\alpha=t_{1}^{i_{1}}\dots t_{n}^{i_{n}}\prod_{{\bold j}}\varepsilon _{{\bold j}
}^{b_{{\bold j}}}\omega^{c}\mod K^{* p^{m}}$,
$b_{{\bold j}},c\in\Bbb Z_p$.

\Item{(2)} 
$
V\left( \{ t_{1},\dots ,t_{n}\} ,\varepsilon
_{{\bold j}}\right) _{m}=1, \quad 
V\left( \{ t_{1},\dots ,t_{n}\} ,\omega\right) _{m}=\zeta _{p^m}$.
\endRoster  

An important element of a Shafarevich basis is
$
\omega(a)=E(as(X))|_{X_n=t_n, \dots , X_1=t_1}$ where
$$E(f(X))=\expp \biggl( \bigl( 1+\frac{\triangle }{p}%
+\frac{\triangle^2 }{p^2}+\cdots \bigr) (f(X))\biggr),$$
$a\in W(K_0)$.

Now take the following elements  as a Shafarevich basis of $K^*/K^{*p^m}$:  

\Roster
\Item{---} elements $t_{1},\dots ,t_{n}$,

\Item{---} elements
$
\varepsilon _{J}=1+\theta t_{n}^{j_{n}}\dots t_{1}^{j_{1}}$ 
where $p\nmid  \text{gcd}\,(j_{1},\dots ,j_{n})$, 
\Item{} $0<(
j_{1},\dots ,j_{n}) <p(e_1,\dots,e_n)/(p-1)$,   
where $(e_1,\dots,e_n)=\bold v(p)$, $\bold v$ is the discrete valuation of rank $n$ associated to $t_{1},\dots ,t_{n}$, 

\Item{---} $\omega=\omega(a)$ where $a$ is an appropriate generator of $W(K_0)/({\bold F}-1)W(K_0)$. 
\endRoster 

Using this basis it  is relatively easy to show that $V(\,\,,\,)_{m}$ is non-degenerate.

In particular, for every $\theta \in \Cal R^*$ there is $\theta'\in\Cal R^*$ such that 
$$V\bigl(\bigl\{1+\theta t_n^{i_n}\dots t_1^{i_1}, t_1,\dots, \widehat{t_l},\dots,t_n\bigr\},
1+\theta' t_n^{pe_n/(p-1)-i_n}\dots t_1^{pe_1/(p-1)-i_1}\bigr)_m=\zeta_{p^m}$$
where $i_l$ is prime to $p$,
 $0<(i_{1},\dots ,i_{n}) <p(e_1,\dots,e_n)/(p-1)$  
and $(e_1,\dots,e_n)=\bold v(p)$. 
\endpf

\th Theorem  2

Every open subgroup $N$ of finite index in $K_{n}^{\tpp}(K)$ such that 
$N\supset p^{m}K_{n}^{\tpp}(K)$ 
is the orthogonal complement \ with respect to $V(\,\,,\,)_m$ of
\ a subgroup in $K^*/K^{*p^m}$.
\endth

\rk Remark

  Given higher local class field theory
one defines the Hilbert symbol
for $l$ such that $l$ is not divisible by $\chr(K)$, 
$\mu_l\le K^*$ as
 $$
(\,\,,\,\,)_{l}\colon K_n(K)/l\times K^*/K^{*\, l}\to \langle \zeta_{l}\rangle,
\quad (x,\beta)_{l}=
\gamma^{\Psi _{K}(x )-1}$$
where $\gamma^l=\beta$, 
$\Psi_{K}\colon  K_n(K) \to \Gal(K^{\ab}/K)$
is the reciprocity map.

If $l$ is prime to $p$, then 
the Hilbert symbol $(\,\,,\,)_{l}$ coincides (up to a sign) with
the $(q-1)/l$th power of the tame symbol of 6.4.2.
If $l=p^m$, then 
the $p^{m}$th Hilbert
symbol
coincides (up to a sign) with 
the symbol $V(\,,\,\,)_{m}$. 
\endrk

 \Bib        References

\rf{A1} V. Abrashkin,
The field of norms functor and the Br\"uckner--Vostokov
formula,
Math. Ann. 308(1997), 5--19.

\rf{A2}  V. Abrashkin, 
   Explicit formulae for the Hilbert symbol of a formal group
over Witt vectors, Izv. Ross. Akad. Nauk Ser. Mat. (1997);
English translation in   Izv. Math. 
   61(1997), 463--515.

\rf{AH1} E. Artin and H. Hasse, 
\"Uber den zweiten Erg\"an\-zun\-g\-s\-satz zum Re\-zi\-p\-ro\-zi\-t\"at\-s\-ge\-setz der
$l$-ten Po\-ten\-z\-res\-te im K\"orper $k_{\zeta}$ der $l$-ten
Einheits\-wur\-zeln und Ober\-k\"orpern von $k_{\zeta}$, 
J. reine angew. Math. 154(1925),  143--148.

\rf{AH2} E. Artin and H. Hasse, 
Die beiden Erg\"anzungssatz zum
Re\-zi\-pr\-zi\-t\"at\-s\-ge\-setz der $l^n$-ten Po\-ten\-z\-re\-s\-te im K\"or\-per
der $l^n$-ten Ein\-heit\-s\-wur\-zeln, 
Abh. Math. Sem. Univ. Hamburg 6(1928),  146--162.

\rf{Be1}   D. Benois, 
   P\'eriodes $p$--adiques et lois de r\'eciprocit\'e explicites, 
   J. reine angew. Math. 
   493(1997),    115--151.

 \rf {Be2} D. Benois, On Iwasawa theory of cristalline representations, Preprint Inst. Experimentelle
 Mathematik (Essen) 1998.

\rf{BK} S. Bloch and K.~Kato, 
   $L$--functions and Tamagawa numbers of motives, 
   \ In \ The Gro\-the\-n\-di\-eck Festschrift, Birkh\"auser
   vol. 1, 1990,     334--400.
  
\rf{Br} H. Br\"uckner,
   Explizites reziprozit\"atsgesetz und Anwendungen, 
   Vorlesungen aus dem Fachbereich Mathematik der Universit\"at Essen,
   1979. 
  
\rf{Cole1} R. F. Coleman, 
   The dilogarithm and the norm residue symbol, \ 
   Bull. Soc. \  France   \   109(1981),    373--402.

\rf{Cole2}  R. F. Coleman,
The arithmetic of Lubin--Tate division towers,
Duke Math. J. 48(1981), 449--466.

\rf {Colm} P. Colmez, Th\'eorie d'Iwasawa des repr\'esentations de de Rham d'un corps local Ann. of Math. 148 (1998), no 2,
485--571.

\rf{D} F. Destrempes, 
   Explicit reciprocity laws for Lubin--Tate modules, 
   J. reine angew. Math., 
   463(1995)   27--47.

\rf{dS} E. de Shalit, 
   The explicit reciprocity law in local class field theory, 
   Duke Math. J.   53(1986),   163--176.

\rf{Fe1}  I. Fesenko, The generalized Hilbert symbol in
the $2$-adic case, 
Vestnik Leningrad. Univ. Mat. Mekh. Astronom. (1985);  
English translation in Vestnik Leningrad Univ. Math. 18(1985), 88--91.

\rf{Fe2} I. Fesenko,
Explicit constructions in local fields,
Thesis, St. Petersburg Univ. 1987.

\rf{Fe3} I. Fesenko,
Generalized Hilbert symbol in multidimensional local fields,
Rings and modules. Limit theorems of probability theory, No. 2
(1988).

\rf{FV} I. Fesenko and S. Vostokov, 
Local Fields and Their Extensions,
AMS, Providence RI, 1993.

\rf{Fo} J.-M. Fontaine,
Appendice: Sur un th\'eor\`eme de Bloch et Kato
(lettre \`a B. Perrin--Riou),
Invent. Math. 115(1994), 151--161.

\rf{Fu} T. Fukaya, 
   Explicit reciprocity laws for $p$--divisible groups
over higher dimensional local fields, 
preprint 1999.

\rf{H} G. Henniart, Sur les lois de r\`eciprocit\'e explicites I,
J. reine angew. Math. 329(1981), 177-203.
  
\rf{I}   K. Iwasawa, 
   On explicit formulas for the norm residue symbols, 
   J. Math. Soc. Japan    20(1968),     151--165.

\rf{Ka1} K.~Kato, 
   The explicit reciprocity law and the cohomology of 
Fontaine--Messing, 
   Bull. Soc. Math. France  
   119(1991),   397--441.

\rf{Ka2} K.~Kato, 
Lectures on the approach to Hasse -Weil  L-functions via
$B_{dR}$, 
Lect. Notes in Math.  1553 (1993), 50--163.

\rf{Ka3} K.~Kato,
Generalized explicit reciprocity laws, Adv. Stud. in Contemporary Math. 
1(1999), 57--126. 

\rf{Kn} M. Kneser,
Zum expliziten Resiprozit\"atsgestz von I. R. Shafarevich,
Math. Nachr. 6(1951), 89--96.

\rf{Ko}  V. Kolyvagin, 
   Formal groups and the norm residue symbol,
Izv. Akad. Nauk Ser. Mat. (1979);
English translation in 
   Math. USSR Izv.  
  15(1980),   289--348.

\rf{Kum}   E. Kummer, 
   \"Uber die allgemeinen Reziprozit\"atsgesetze der Potenzreste, 
    J.  reine  angew.  Math. 
   56(1858),     270--279.

\rf{Kur1} M. Kurihara, Appendix: Computation of the syntomic regulator
in the cyclotomic case,  
Invent. Math. 99 (1990), 313--320.

\rf{Kur2}  M. Kurihara, 
   The exponential homomorphism for the Milnor $K$--groups
and an explicit reciprocity law,
J. reine angew. Math. 498(1998), 201--221.

\rf{PR} B. Perrin-Riou, 
   Theorie d'Iwasawa des representations  $p$--adiques
sur un corps local, 
   Invent.  Math.   115(1994),   81--149.

\rf{Se}  Sh. Sen, 
   On explicit reciprocity laws I; II, 
   J. reine angew. Math.  
   313(1980),  1--26; 323(1981), 68--87.

\rf{Sh} I. R. Shafarevich, 
A general reciprocity law, 
Mat. Sb. (1950); 
English translation in Amer. Math. Soc. Transl. Ser. 2, 4(1956),  73--106.

\rf{Su} Y. Sueyoshi, \ 
Explicit reciprocity laws on relative Lubin--Tate groups, \ 
Acta Arithm. \ 55(1990), 291--299.

\rf{V1} S.  V.  Vostokov,
   An explicit form of the reciprocity law,
 Izv. Akad. Nauk SSSR Ser. Mat. (1978); English translation in 
   Math.  USSR Izv.   13(1979), 557--588.
    
\rf{V2} S.  V.  Vostokov,
   A norm pairing in formal modules, 
 Izv. Akad. Nauk SSSR Ser. Mat. (1979);
English translation in   Math.  USSR Izv.  15(1980),   25--52.

\rf{V3}  S. V. Vostokov,
   Symbols on formal groups, 
Izv. Akad. Nauk SSSR Ser. Mat. (1981);
English translation in 
  Math.  USSR Izv.      19(1982),   261--284.

\rf{V4} 
  S. V. Vostokov,  The Hilbert symbol for Lubin--Tate formal groups I, 
    Zap.  Nauchn.  Sem.  Leningrad.  Otdel.  Mat.  Inst.  Steklov.
(LOMI) (1982); 
English translation in  J.  Soviet Math.     27(1984), 2885--2901.

\rf{V5}  S. V. Vostokov, 
   Explicit construction of class field theory for a multidimensional
local field, 
   Izv.  Akad.  Nauk SSSR Ser.  Mat.   (1985) no.2;
English translation in  Math.  USSR Izv.      26(1986),  263--288.

\rf{V6}   S. V.  Vostokov, 
   The pairing on $K$-groups in fields of valuation
of rank $n$, 
   Trudy Sankt-Peterb. Mat. Obschestva (1995); 
English translation in   Amer. Math. Soc. Transl. Ser. 2 
   165(1995),    111--148.

\rf{V7} S. V. Vostokov, 
   Hilbert pairing on a multidimensional complete field
   Trudy Mat. Inst. Stek\-lo\-va, 
  (1995); 
English translation in  Proc. Steklov Inst. of Math.
  208(1995), 
  72--83.

\rf{VB} S.  V.  Vostokov and  D.  G.  Benois, 
   Norm pairing in formal groups and Galois representations, 
    Algebra i Analiz (1990);
English translation in  Leningrad Math.  J.     2(1991), 1221--1249. 
  
\rf{VD1}   S. V. Vostokov and O. V. Demchenko, 
   Explicit form of Hilbert
pairing for relative Lubin-Tate formal groups, 
   Zap. Nauch. Sem. POMI 
   (1995); 
English translation in  J. Math. Sci.  Ser. 2
   89(1998),    1105--1107.

\rf{VD2} S.V.Vostokov and O.V. Demchenko, Explicit formula of the Hilbert symbol for
Honda formal group, Zap. Nauch. Sem. POMI 272(2000) (in Russian).

\rf{VF}  S.  V.  Vostokov and I. B.  Fesenko, 
   The Hilbert symbol for Lubin--Tate formal groups II, 
    Zap.  Nauchn.  Sem.  Leningrad.  Otdel.  Mat.  Inst.  Steklov (LOMI), 
  (1983); 
English translation in  J.  Soviet Math.     30(1985),    1854--1862.

\rf{W}    A. Wiles, 
   Higher explicit reciprocity laws, 
   Ann. Math. 
   107(1978),   235--254.
\endBib

\Coordinates

Department of Mathematics \ 
St. Petersburg University

Bibliotechnaya pl. 2\ 
Staryj Petergof,
 198904 St. Petersburg \ Russia

E-mail: sergei\@vostokov.usr.pu.ru

\endCoordinates
\vfill
\pagebreak

\end

%% file: m3-macs.tex
\expandafter\ifx\csname mthreemacsloaded\endcsname\relax\else \fi

\magnification1100
\input amstex


 \catcode`\@=11
 \let\wlog@ld\wlog
 \def\wlog#1{\relax}

 \newif\ifIN@
 \def\m@rker{\m@@rker}
 \def\IN@{\expandafter\INN@\expandafter}
 \long\def\INN@0#1@#2@{\long\def\NI@##1#1##2##3\ENDNI@
    {\ifx\m@rker##2\IN@false\else\IN@true\fi}%
     \expandafter\NI@#2@@#1\m@rker\ENDNI@}
  \newtoks\Initialtoks@  \newtoks\Terminaltoks@
  \def\SPLIT@{\expandafter\SPLITT@\expandafter}
  \def\SPLITT@0#1@#2@{\def\TTILPS@##1#1##2@{%
     \Initialtoks@{##1}\Terminaltoks@{##2}}\expandafter\TTILPS@#2@}
  \newtoks\Trimtoks@

 \def\ForeTrim@{\expandafter\ForeTrim@@\expandafter}
 \def\ForePrim@0 #1@{\Trimtoks@{#1}}
 \def\ForeTrim@@0#1@{\IN@0\m@rker. @\m@rker.#1@%
     \ifIN@\ForePrim@0#1@%
     \else\Trimtoks@\expandafter{#1}\fi}
 
  \def\Trim@0#1@{%
      \ForeTrim@0#1@%
      \IN@0 @\the\Trimtoks@ @%
        \ifIN@
             \SPLIT@0 @\the\Trimtoks@ @\Trimtoks@\Initialtoks@
             \IN@0\the\Terminaltoks@ @ @%
                 \ifIN@
                 \else \Trimtoks@ {FigNameWithSpace}%
                 \fi
        \fi
      }

  \font\titlebold=cmbx12 scaled 1200
  \font\twelvebold=cmbx12
  \font\tenbold=cmbx10
  \font\ninebold=cmbx9
  \font\sevenbold=cmbx7
  \font\fivebold=cmbx5

  \input amssym.def \input amssym
     \font\titlemsa=msam10 at 14.4pt
     \font\titlemsb=msbm10 at 14.4pt
     \font\titleeufm=eufm10 at 14.4pt
     \font\twelvemsa=msam10 scaled 1200
     \font\twelvemsb=msbm10 scaled 1200
     \font\twelveeufm=eufm10 scaled 1200
     \font\ninemsa=msam9
     \font\ninemsb=msbm9
     \font\nineeufm=eufm9

   \ifx\cyrfam\undefined
   \else
     \immediate\write16{}%
     \message{ !!! cyr fonts already defined. !!! }
     \message{ --- edit out superfluous font defs? }
   \fi
   \newfam\cyrfam
       \font\titlecyr=wncyr10 scaled 1440 
       \font\twelvecyr=wncyr10 scaled 1200
       \font\tencyr=wncyr10
       \font\ninecyr=wncyr9
       \font\sevencyr=wncyr7
       \font\sixcyr=wncyr6

   \newfam\eusmfam
       \font\titleeusm=eusm10 scaled 1440
       \font\twelveeusm=eusm10 scaled 1200
       \font\teneusm=eusm10
       \font\nineeusm=eusm9
       \font\seveneusm=eusm7
       
       \font\fiveeusm=eusm5

\let\Cal\cal

    \font\ninemrm=cmr9 
    \font\ninei=cmmi9
    \font\ninesy=cmsy9 
    \skewchar\ninei='177
    \skewchar\ninesy='60

  \font\twelvemrm=cmr10 at 12pt 
  \font\twelvei=cmmi10 at 12pt
  \font\twelvesy=cmsy10 at 12pt

  \font\titlemrm=cmr10 at 14.4pt 
  \font\titlei=cmmi10 at 14.4pt
  \font\titlesy=cmsy10 at 14.4pt


  \def\Smallfonts{\ninepoint}

  \def\Hfont{\titlepoint\bf}
  \def\Authorfont{\twelvepoint\it}
  \def\HHfont{\twelvepoint\bf}
  \def\HHHfont{\bf}
  \def\Bibfont{\tenbf}
  \def\Coordfont{\nineit }

  \def \thfont {\bf }
  \def \pffont {\it\itSpacing }
  \def \rkfont {\bf }
  \def \dffont {\bf }
  \def \egfont {\bf }

 \def\ninepoint{%
  \def\rm{\fam0\ninerm}%
    \textfont0=\ninemrm  \scriptfont0=\sevenrm  \scriptscriptfont0=\fiverm
    \textfont1=\ninei    \scriptfont1=\seveni   \scriptscriptfont1=\fivei
  \def\mit{\fam1\ninei}%
  \def\oldstyle{\fam1\ninei}%
    \textfont2=\ninesy   \scriptfont2=\sevensy  \scriptscriptfont2=\fivesy
    \textfont3=\tenex    \scriptfont3=\tenex    \scriptscriptfont3=\tenex
  \def\it{\fam\itfam\nineit}%
    \textfont\itfam=\nineit
  \def\bf{\ifmmode\fam\bffam\else\ninebf\fi}%
    \textfont\bffam=\ninebold 
    \scriptfont\bffam=\sevenbold 
    \scriptscriptfont\bffam=\fivebold%
  \def\msa{\fam\msafam\ninemsa}%
    \textfont\msafam=\ninemsa 
    \scriptfont\msafam=\sevenmsa
    \scriptscriptfont\msafam=\fivemsa%
  \def\msb{\fam\msbfam\ninemsb}%
    \textfont\msbfam=\ninemsb%
    \scriptfont\msbfam=\sevenmsb%
    \scriptscriptfont\msbfam=\fivemsb%
  \def\eufm{\fam\eufmfam\nineeufm}%
    \textfont\eufmfam=\nineeufm
    \scriptfont\eufmfam=\seveneufm
    \scriptscriptfont\eufmfam=\fiveeufm
   \def\eusm{\fam\eusmfam\nineeusm}%
     \textfont\eusmfam=\nineeusm
     \scriptfont\eusmfam=\seveneusm
     \scriptscriptfont\eusmfam=\fiveeusm
   \def\cyr{\fam\cyrfam\ninecyr}%
     \textfont\cyrfam=\ninecyr
     \scriptfont\cyrfam=\sevencyr
     \scriptscriptfont\cyrfam=\sixcyr
  \setbox\strutbox=\hbox{\vrule
      height7pt depth3pt width0pt}%
   \baselineskip=10.8pt\rm}

 \let\eightpoint\ninepoint 

 \def\tenpoint{%
  \def\rm{\fam0\tenrm}%
    \textfont0=\tenmrm \scriptfont0=\sevenrm \scriptscriptfont0=\fiverm%
  \def\mit{\fam1\teni}%
  \def\oldstyle{\fam1\teni}%
    \textfont1=\teni   \scriptfont1=\seveni  \scriptscriptfont1=\fivei%
    \textfont2=\tensy  \scriptfont2=\sevensy \scriptscriptfont2=\fivesy%
    \textfont3=\tenex  \scriptfont3=\tenex   \scriptscriptfont3=\tenex%
  \def\it{\fam\itfam\tenit}%
    \textfont\itfam=\tenit%
  \def\bf{\ifmmode\fam\bffam\else\tenbf\fi}%
    \textfont\bffam=\tenbold
    \scriptfont\bffam=\sevenbold%
    \scriptscriptfont\bffam=\fivebold%
  \def\msa{\fam\msafam\tenmsa}%
    \textfont\msafam=\tenmsa%
    \scriptfont\msafam=\sevenmsa%
    \scriptscriptfont\msafam=\fivemsa%
  \def\msb{\fam\msbfam\tenmsb}%
    \textfont\msbfam=\tenmsb%
    \scriptfont\msbfam=\sevenmsb%
    \scriptscriptfont\msbfam=\fivemsb%
  \def\eufm{\fam\eufmfam\teneufm}%
   \textfont\eufmfam=\teneufm
   \scriptfont\eufmfam=\seveneufm
   \scriptscriptfont\eufmfam=\fiveeufm
   \def\eusm{\fam\eusmfam\teneusm}%
    \textfont\eusmfam=\teneusm
    \scriptfont\eusmfam=\seveneusm
    \scriptscriptfont\eusmfam=\fiveeusm
   \def\cyr{\fam\cyrfam\tencyr}%
    \textfont\cyrfam=\tencyr
    \scriptfont\cyrfam=\sevencyr
    \scriptscriptfont\cyrfam=\sixcyr
  \setbox\strutbox=\hbox{\vrule %
      height8.5pt depth3.5ptwidth0pt}%
  \baselineskip=\StdBaselineskip\rm}

 \def\twelvepoint{%
  \def\rm{\fam0\twelverm}%
    \textfont0=\twelvemrm \scriptfont0=\tenmrm \scriptscriptfont0=\sevenrm
    \textfont1=\twelvei   \scriptfont1=\teni   \scriptscriptfont1=\seveni
  \def\mit{\fam1\twelvei}%
  \def\oldstyle{\fam1\twelvei}%
    \textfont2=\twelvesy  \scriptfont2=\tensy  \scriptscriptfont2=\sevensy
    \textfont3=\tenex  \scriptfont3=\tenex  \scriptscriptfont3=\tenex
  \def\it{\fam\itfam\twelveit}%
    \textfont\itfam=\twelveit
  \def\bf{\ifmmode\fam\bffam\else\twelvebf\fi}%
    \textfont\bffam=\twelvebold
    \scriptfont\bffam=\tenbold%
    \scriptscriptfont\bffam=\sevenbold%
  \def\msa{\fam\msafam\twelvemsa}%
    \textfont\msafam=\twelvemsa%
    \scriptfont\msafam=\tenmsa%
    \scriptscriptfont\msafam=\sevenmsa%
  \def\msb{\fam\msbfam\twelvemsb}%
    \textfont\msbfam=\twelvemsb%
    \scriptfont\msbfam=\tenmsb%
    \scriptscriptfont\msbfam=\sevenmsb%
  \def\eufm{\fam\eufmfam\twelveeufm}%
   \textfont\eufmfam=\twelveeufm
   \scriptfont\eufmfam=\teneufm
   \scriptscriptfont\eufmfam=\seveneufm
   \def\eusm{\fam\eusmfam\twelveeusm}%
    \textfont\eusmfam=\twelveeusm
    \scriptfont\eusmfam=\teneusm
    \scriptscriptfont\eusmfam=\seveneusm
   \def\cyr{\fam\cyrfam\tencyr}%
    \textfont\cyrfam=\twelvecyr
    \scriptfont\cyrfam=\tencyr
    \scriptscriptfont\cyrfam=\sevencyr
  \setbox\strutbox=\hbox{\vrule
      height10.2pt depth4.55pt width0pt}%
  \baselineskip=14pt\rm}

 \def\titlepoint{%
    \textfont0=\titlemrm \scriptfont0=\twelvemrm \scriptscriptfont0=\tenmrm
    \textfont1=\titlei   \scriptfont1=\twelvei   \scriptscriptfont1=\teni
  \def\mit{\fam1\titlei}%
  \def\oldstyle{\fam1\titlei}%
    \textfont2=\titlesy  \scriptfont2=\twelvesy  \scriptscriptfont2=\tensy
    \textfont3=\tenex
    \scriptfont3=\tenex
    \scriptscriptfont3=\tenex
  \def\it{\fam\itfam\titleit}%
    \textfont\itfam=\titleit
  \def\bf{\ifmmode\fam\bffam\else\titlebf\fi}%
    \textfont\bffam=\titlebold
    \scriptfont\bffam=\twelvebold%
    \scriptscriptfont\bffam=\tenbold%
  \def\msa{\fam\msafam\titlemsa}%
    \textfont\msafam=\titlemsa%
    \scriptfont\msafam=\twelvemsa%
    \scriptscriptfont\msafam=\tenmsa%
  \def\msb{\fam\msbfam\titlemsb}%
    \textfont\msbfam=\titlemsb%
    \scriptfont\msbfam=\twelvemsb%
    \scriptscriptfont\msbfam=\tenmsb%
  \def\eufm{\fam\eufmfam\titleeufm}%
    \textfont\eufmfam=\titleeufm
    \scriptfont\eufmfam=\twelveeufm
    \scriptscriptfont\eufmfam=\teneufm
   \def\eusm{\fam\eusmfam\titleeusm}%
     \textfont\eusmfam=\titleeusm
     \scriptfont\eusmfam=\twelveeusm
     \scriptscriptfont\eusmfam=\teneusm
   \def\cyr{\fam\cyrfam\tencyr}%
    \textfont\cyrfam=\titlecyr
    \scriptfont\cyrfam=\twelvecyr
    \scriptscriptfont\cyrfam=\tencyr
  \setbox\strutbox=\hbox{\vrule
      height12.3pt depth5.54pt width0pt}%
  \baselineskip=16pt\rm}

\newbox\AuthorBox\newbox\TitleBox
\newbox\TFLinebox
\newbox\FLinebox
\newbox\HLinebox
\def\SetTFLinebox#1{\setbox\TFLinebox=\hbox{#1}}
\def\SetFLinebox#1{\setbox\FLinebox=\hbox{#1}}
\def\SetHLinebox#1{\setbox\HLinebox=\hbox{#1}}

 \def\SetAuthorHead#1{%
     \setbox\AuthorBox=\hbox{\ninepoint \it 
           \ignorespaces\frenchspacing#1\unskip}}
 \def\SetTitleHead#1{%
     \setbox\TitleBox=\hbox{\ninepoint \it
           \ignorespaces\frenchspacing#1\unskip}}

  \def\itSpacing{\relax}
  \def\itSpacingOff{\relax}


 \def\Hrule{\hrule width0pt height0pt}

  \newskip\ProcSkip \ProcSkip 8pt plus2pt minus2pt

 \newskip\LastSkip
 \def\SaveLastSkip{\LastSkip\lastskip}
 \def\RestoreLastSkip{\vskip-\LastSkip\vskip\LastSkip}

 \def\NoindentAfter{\everypar={\setbox0=\lastbox\everypar={}}}

 \long\def\H#1\par#2\par{\notenumber=0 \titlepagetrue%
    {
    \baselineskip=20pt
    \parindent=0pt\parskip=0pt\frenchspacing
    \leftskip=0pt plus .2\hsize minus .3\hsize
    \rightskip=0pt plus .2\hsize minus .3\hsize
 \def\\{\unskip\break}%
    \pretolerance=10000 \Hfont #1\unskip\break
     \vskip7pt\Hrule
\hfill \Authorfont #2\hfill\hfill\unskip}
    \vskip48pt plus 4pt minus 4pt
    \par\NoindentAfter\rm}

 \long\def\Hi#1\par#2\par{\notenumber=0 \titlepagetrue%
    {  \baselineskip=0pt  \parindent=0pt\parskip=0pt\frenchspacing
    \leftskip=0pt plus .2\hsize minus .3\hsize
    \rightskip=0pt plus .2\hsize minus .3\hsize
}
    \rm}


 \newdimen\PageRemainder
  \def\SetPageRemainder{
     \PageRemainder=\pagegoal
     \ifdim\PageRemainder=\maxdimen\PageRemainder=\vsize
     \else\advance\PageRemainder by -1\pagetotal\fi}

  \def\Rpt@{}\def\Rpt@@{}

  \long\def\HH#1\par{\par
  \SaveLastSkip\removelastskip\goodbreak
  \ifdim\LastSkip<30pt 
     \LastSkip 30pt
plus 3pt minus 2pt\fi
  \SetPageRemainder\advance\PageRemainder-\LastSkip
  \ifdim\PageRemainder<150pt
       \edef\Rpt@{remain = \the\PageRemainder\noexpand\\
                pagetotal=\the\pagetotal\noexpand\\
                           pagegoal=\the\pagegoal}%
          \fi
   \ifdim\PageRemainder<65pt 
       \ifdim\PageRemainder > 0pt
          \edef\Rpt@@{\noexpand\\
                      Had HH PageRemainder$<$\relax 65pt\noexpand\\
                      Hence forced break!}%
     \vskip 0pt plus .2\PageRemainder\eject 
    \fi\fi
    \vskip\LastSkip\Hrule 
    \pretolerance=10000\rightskip=0pt plus 3em
    \hangafter1 \hangindent=2.2em%
    \noindent
    \HHfont \unskip \Ednote{\Rpt@\Rpt@@}%
            \def\Rpt@{}\def\Rpt@@{}%
            \ignorespaces
            #1\par\rightskip=0pt\pretolerance=\StdPretolerance%
    \NoindentAfter
\tenpoint\rm%
     \medskip \vskip\ProcSkip}

  \long\def\HHH#1\par{\par%
  \SaveLastSkip\removelastskip\goodbreak
  \ifdim\LastSkip<\ProcSkip%
     \LastSkip\ProcSkip\fi
  \SetPageRemainder\advance\PageRemainder-\LastSkip
  \ifdim\PageRemainder<150pt
       \edef\Rpt@{remain = \the\PageRemainder\noexpand\\
                pagetotal=\the\pagetotal\noexpand\\
                           pagegoal=\the\pagegoal}%
       \fi
   \ifdim\PageRemainder<48pt  
        \ifdim\PageRemainder > 0pt
             \edef\Rpt@@{\noexpand\\
                      Had HHH PageRemainder$<$\relax48pt\noexpand\\
                      Hence forced break!}%
       \vskip 0pt plus .2\PageRemainder\eject 
      \fi\fi
   \vskip\LastSkip\par\noindent
   \HHHfont \unskip\Ednote{\Rpt@\Rpt@@}%
  \def\Rpt@{}\def\Rpt@@{}%
  \ignorespaces
   #1\unskip.\quad\rm\ignorespaces
   \ignorepars}

  \long\def\ignorepars#1\par{\def\Test{#1}%
     \ifx\Test\Empty\def\This{\ignorepars}%
        \else\def\This{\Test\par}\fi
           \This}
  \def\Empty{}

 \def\Abstract#1\par{\bgroup\Smallfonts\narrower\HHH #1\par}
 \def\endAbstract{\par\egroup}


 \def\ProcBreak{\par%
    \ifdim\lastskip<8pt%
    \removelastskip%
    \penalty-200\vskip\ProcSkip\fi}

 \def\th#1\par{\ProcBreak \noindent
   {\thfont\ignorespaces
    #1\unskip.}\it\itSpacing\kern.4em\ignorepars}

 \def\endth{\ProcBreak\rm\itSpacingOff }


 \def\pf#1\par{\ProcBreak %
    \noindent\pffont#1\unskip.\rm\itSpacingOff{\kern .7em}\ignorepars}

 \def\endpf{\medskip \ProcBreak } 

  \def\qedbox{\hbox{\vbox{
    \hrule width0.2cm height0.2pt
    \hbox to 0.2cm{\vrule height 0.2cm width 0.2pt
             \hfil\vrule height0.2cm width 0.2pt}
    \hrule width0.2cm height 0.2pt}\kern1pt}}

  \def\qed{\ifmmode\qedbox
    \else\unskip\ \hglue0mm\hfill\qedbox\ProcBreak\fi}

  \def \rk #1\par{\ProcBreak
     \noindent{\rkfont\ignorespaces #1\unskip.}%
     \rm\kern.6em\ignorepars}

  \def \endrk {\medskip\ProcBreak }

  \def \df #1\par{\ProcBreak
     \noindent{\dffont\unskip\ignorespaces #1\unskip.}%
     \rm\kern.6em\ignorepars}

  \def \enddf {\medskip\ProcBreak }

  \def \eg #1\par{\ProcBreak
     \noindent\egfont\unskip\ignorespaces #1\unskip.
     \rm\kern.6em\ignorepars}

  \newdimen\Overhang

   \def\MaxTag@#1#2#3#4#5{\setbox0=\hbox{#4\ignorespaces#2\unskip}%
     \dimen0=\wd0\advance\dimen0 by#3
     \ifdim\dimen0<#5\relax\dimen0=#5\fi
     \expandafter\edef\csname #1Hang\endcsname{\the\dimen0}}

 \def\MaxItemTag#1{\MaxTag@{Item}{#1}{.4em}{\ItemStyle}{\parindent}}%
 \def\MaxItemItemTag#1{%
        \MaxTag@{ItemItem}{#1}{.4em}{\ItemItemStyle}{\parindent}}
 \def\MaxNrTag#1{\MaxTag@{Nr}{#1}{.5em}{\NrStyle}{\parindent}}
 \def\MaxReferenceTag#1{%
        \MaxTag@{Reference}{[#1]}{.6em}{\ninerm}{\parindent}}
 \def\MaxFootTag#1{\MaxTag@{Foot}{#1}{.4em}{\ninerm}{\z@}}

  \def\SetOverhang@{\Overhang=.8\dimen0%
     \advance\Overhang by \wd0\relax
     \ifdim\Overhang>\hangindent\relax
       \advance\Overhang by .25\dimen0%
       \Ednote{Tag is pushing text.}\osumess{Tag is pushing text.}%
     \else\Overhang=\hangindent
     \fi}

   \def\Item#1{\par\noindent
      \hangafter1\hangindent=\ItemHang
      \setbox0=\hbox{\ItemStyle\ignorespaces#1\unskip}%
      \dimen0=.4em\SetOverhang@
      \rlap{\box0}\kern\Overhang\ignorespaces}

   \def\ItemItem#1{\par\noindent
      \hangafter1\hangindent=\ItemItemHang
      \setbox0=\hbox{\ItemItemStyle\ignorespaces#1\unskip}%
      \dimen0=.4em\SetOverhang@
      \advance\hangindent by \ItemHang
      \kern\ItemHang\rlap{\box0}%
      \kern\Overhang\ignorespaces}

  \def\Nr#1{\par\noindent\hangindent=\NrHang 
    \setbox0=\hbox{\NrStyle\ignorespaces#1\unskip}%
    \dimen0=.5em\SetOverhang@
    \rlap{\box0}\kern\Overhang
    \hangindent=\z@\ignorespaces}

   \newskip\Rosterskip\Rosterskip 1pt plus1pt 
   \def\Roster{\par\ifdim\lastskip<\Rosterskip\removelastskip\vskip\Rosterskip\fi
    \bgroup}
   \def\endRoster{\par\global\edef\LastSkip@{\the\lastskip}\removelastskip
       \egroup\penalty-50\LastSkip\LastSkip@\relax
       \ifdim\LastSkip<\Rosterskip\LastSkip\Rosterskip\fi
       \vskip\LastSkip}




 \def\cite#1{
    \def\nextiii@##1,##2\end@{{\frenchspacing\rm 
      \lBr\ignorespaces##1\unskip{\rm,~\ignorespaces##2}\rBr}}%
    \IN@0,@#1@%
    \ifIN@\def\next{\nextiii@#1\end@}\else
    \def\next{{\rm\lBr#1\rBr}}\fi\next}


   \def \Bib#1\par{%
       \par\removelastskip\SetPageRemainder
       \ifdim\PageRemainder < 97pt
        \ifdim\PageRemainder > 0pt
        \vfill\eject
       \fi\fi
    \ProcBreak \par\begingroup\parskip=0 pt%
    \goodbreak \vskip 15 pt plus 10 pt
    \noindent\null\hfill\Bibfont
      \ignorespaces #1\unskip\hfill\null\par 
    \frenchspacing \Smallfonts\rm
    \parskip=2.5 pt plus 1 pt minus.5pt%
    \nobreak\vskip 12pt plus 2pt minus2pt\nobreak
    \leftskip=0 pt \baselineskip=10.5pt}

 \def\ReferenceTagSlide{0em}
  \def\ReferenceTagGap{.5em}

  \def \rf#1{\par\noindent
     \hangafter1\hangindent=\ReferenceHang      
     \setbox0=\hbox{\ninerm[\ignorespaces#1\unskip]}%
     \dimen0=\ReferenceTagGap\SetOverhang@
     \rlap{\kern\ReferenceTagSlide\box0}%
     \kern\Overhang\ignorespaces}

  \def\ref#1\par#2\par#3\par#4\par{%
     \rf{#1}#2\unskip,\ #3\unskip,\
     #4\unskip.}

  \def\endBib{\par\endgroup\vskip 12pt minus 6pt }


  \long\def\Coordinates#1\endCoordinates{
 {\par\vskip4pt\def\\{\unskip, }\Coordfont\baselineskip10.5pt\noindent#1}}

 \def\pagecontents{
  \gdef\Pagetot@l{\pagetotal}
  \ifvoid\TRMargIns\else
    \rlap{\kern\hsize\kern10pt\vbox to 0pt{%
         \box\TRMargIns\vss}}\fi
  \ifvoid\topins\else\unvbox\topins\fi
   \dimen@=\dp\@cclv \unvbox\@cclv 
   \ifvoid\footins\else 
     \vskip\skip\footins
     \footnoterule
     \unvbox\footins\fi
   \ifr@ggedbottom \kern-\dimen@ \vfil \fi}


 \newcount\Ht 

 \def \Acc{\expandafter } 

 \def\swthat{\raise -1.1 ex\hbox{\sam$\widehat{}$}}
 \def\swttilde{\raise -1.2 ex\hbox{\sam$\widetilde{}$}}
 \def \overdot{{\raise .2 ex \hbox to 0pt {\hss\bf\smash{.}\hss}}}
 \def \overcircle{{\raise .1 ex \hbox to 0pt
    {\sam$\eightpoint\scriptstyle\hss\circ\hss$}}}

 \def \Mathaccent#1#2{{\sam 
  \setbox4=\hbox{$\vphantom{#2}$}
  \Ht=\ht4 
  \setbox5=\hbox{${#1}$}
  \setbox6=\hbox{${#2}$}
  \setbox7=\hbox to .5\wd6{}
  \copy7\kern .1\Ht \raise\Ht sp\hbox{\copy5}\kern-.1\Ht
  \copy7\llap{\box6}
  }}

  \def\SwtCheck #1{
        \ifmmode \check{#1}%
                \else \v {#1}%
                \fi}

 \def\barpartial {%
   \kern .17 em
    \overline {\kern -.17 em\partial\kern-.03 em}%
    \kern .03 em}

 
  \def\Overline#1{\setbox1=\hbox{\sam ${#1}$}%
      \ifdim \wd1 > 6pt
    \kern .11 em
    \overline {\kern -.11 em#1\kern-.14 em}
    \kern .14 em
  \else
    \kern .03 em
    \overline {\kern -.03 em#1\kern-.04 em}
    \kern .04 em
  \fi}

 \def\SOverline#1{\setbox1=\hbox{\sam ${#1}$}%
      \ifdim \wd1 > 7pt
    \kern .22 em
    \overline {\kern -.22 em#1\kern-.09 em}%
    \kern .09 em
  \else
    \kern .10 em
    \overline {\kern -.10 em#1\kern-.04 em}%
    \kern .04 em
  \fi}


 \def\Underline#1{\setbox1=\hbox{\sam ${#1}$}%
      \ifdim \wd1 > 6pt
    \kern .11 em
    \underline {\kern -.11 em#1\kern-.14 em}
    \kern .14 em
  \else
    \kern .03 em
    \underline {\kern -.03 em#1\kern-.04 em}
    \kern .04 em
  \fi}

 \def\SUnderline#1{\setbox1=\hbox{\sam ${#1}$}%
      \ifdim \wd1 > 7pt
    \kern .04 em
    \underline {\kern -.04 em#1\kern-.2 em}%
    \kern .2 em
  \else
    \kern .0 em
    \underline {\kern -.0 em#1\kern-.15 em}%
    \kern .15 em
  \fi}


 \def \Blackbox
   {\leavevmode\hskip .3pt \vbox
   {\hrule height 5pt\hbox{\hskip 4.5pt}}\hskip .5pt}

 \def \XX{\Blackbox\kern.5pt\Blackbox} 

  \def\.{.\kern1pt}

    \def\Hyphen{\edef\this{\the\hyphenchar\font}%
          \hyphenchar\font=-1\char\this\hyphenchar\font=\this}

 \ifx\undefined\text
  \def\text#1{\hbox{\rm #1}}\fi 



   \everymath{}  

  \def\PassMath@@{\aftergroup\AfterMath@} 

 \let\PassMath@\PassMath@@

 \def\AfterMath@{\futurelet\next\AfterMathMole@}

 \def\AfterMathMole@{
      \ifcat\next\space
          \def\this{}
      \else
      \ifcat\next\egroup %
        \def\this{\osumess{Handset mathsurround?? ---(see dollar brace)}}%
      \else
      \def\this{\AAfterMath@}
      \fi\fi
      \this}

 \def\hyphen@{-}
 \def\paren@{)}
 \def\apostr@{'}

 \def\MSC#1{\kern-.8\mathsurround#1\kern.8\mathsurround}

 \def\AAfterMath@#1{\def\Next{#1}
    \IN@0\Next @,.;:!?\relax @%
    \ifIN@\def\this{\MSC{\Next}}%
    \else
    \ifx\Next\hyphen@\def\this{\futurelet\next\AfterHyphen@}%
    \else
    \ifx\Next\paren@\def\this{#1}%
    \else 
    \ifx\Next\apostr@\def\this{#1}%
    \else \def\this{\osumess{Handset mathsurround??}%
                 #1}\fi\fi\fi\fi
    \this}

 \def\AfterHyphen@#1{\def\Next{#1}%
   \ifx\Next\hyphen@\def\this{--}\else
   \ifcat\next\space%
   \def\this{\kern-\mathsurround\kern.05em- \Next}\else
   \def\this{\kern-\mathsurround\kern.05em\Hyphen\Next}\fi\fi\this}

 \def\sam{\mathsurround=\z@\let\PassMath@\relax}  %
 \def\mas{\mathsurround=\StdMathsurround\let\PassMath@\PassMath@@}
 
 \def\Mas{\mathsurround=\StdMathsurround
                \everymath{\PassMath@}\let\PassMath@\PassMath@@}

 \def\m@th{\mathsurround=\z@\everymath{}}

 \def\m@@th{\mathsurround=\z@\everymath={}\let\m@th\relax}

\def\underbar#1{$\setbox\z@\hbox{#1}\dp\z@\z@
      \m@th \underline{\box\z@}$\relax}

\def\mathhexbox#1#2#3{\leavevmode
  \hbox{\m@@th$\m@th \mathchar"#1#2#3$}}

\def\dots{\relax\ifmmode\ldots\else$\m@th\ldots\,$\relax\fi}

\def\dotfill{\cleaders\hbox{\m@@th$\m@th \mkern1.5mu.\mkern1.5mu$}\hfill}
\def\rightarrowfill{$\m@th\mathord-\mkern-6mu%
  \cleaders\hbox{\m@@th$\mkern-2mu\mathord-\mkern-2mu$}\hfill
  \mkern-6mu\mathord\rightarrow$\relax}
\def\leftarrowfill{$\m@th\mathord\leftarrow\mkern-6mu%
  \cleaders\hbox{\m@@th$\mkern-2mu\mathord-\mkern-2mu$}\hfill
  \mkern-6mu\mathord-$\relax}

\def\downbracefill{$\m@th\braceld\leaders\vrule\hfill\braceru
  \bracelu\leaders\vrule\hfill\bracerd$\relax}
\def\upbracefill{$\m@th\bracelu\leaders\vrule\hfill\bracerd
  \braceld\leaders\vrule\hfill\braceru$\relax}

\def\angle{{\vbox{\m@@th\ialign{$\m@th\scriptstyle##$\crcr
      \not\mathrel{\mkern14mu}\crcr
      \noalign{\nointerlineskip}
      \mkern2.5mu\leaders\hrule height.34pt\hfill\mkern2.5mu\crcr}}}}

\def\big#1{{\m@@th\hbox{$\left#1\vbox to8.5\p@{}\right.\n@space$}}}
\def\Big#1{{\m@@th\hbox{$\left#1\vbox to11.5\p@{}\right.\n@space$}}}
\def\bigg#1{{\m@@th\hbox{$\left#1\vbox to14.5\p@{}\right.\n@space$}}}
\def\Bigg#1{{\m@@th\hbox{$\left#1\vbox to17.5\p@{}\right.\n@space$}}}
\def\n@space{\nulldelimiterspace\z@ \m@th}

\def\root#1\of{\setbox\rootbox\hbox{\m@@th$\m@th\scriptscriptstyle{#1}$}
  \mathpalette\r@@t}
\def\r@@t#1#2{\setbox\z@\hbox{\m@@th$\m@th#1\sqrt{#2}$\relax}
  \dimen@\ht\z@ \advance\dimen@-\dp\z@
  \mkern5mu\raise.6\dimen@\copy\rootbox \mkern-10mu \box\z@}

\def\mathph@nt#1#2{\setbox\z@\hbox{\m@@th$\m@th#1{#2}$}\finph@nt}

\def\mathsm@sh#1#2{\setbox\z@\hbox{\m@@th$\m@th#1{#2}$}\finsm@sh}

\def\@vereq#1#2{\lower.5\p@\vbox{\m@@th\baselineskip\z@skip\lineskip-.5\p@
    \ialign{$\m@th#1\hfil##\hfil$\crcr#2\crcr=\crcr}}}

\def\mathpalette#1#2{\sam\mathchoice{#1\displaystyle{#2}}%
  {#1\textstyle{#2}}{#1\scriptstyle{#2}}{#1\scriptscriptstyle{#2}}\mas}

\def\widehat#1{\setbox\z@\hbox{\sam$#1$}%
 \ifdim\wd\z@>\tw@ em\mathaccent"0\msbfam@5B{#1}%
 \else\mathaccent"0362{#1}\fi}
\def\widetilde#1{\setbox\z@\hbox{\sam$#1$}%
 \ifdim\wd\z@>\tw@ em\mathaccent"0\msbfam@5D{#1}%
 \else\mathaccent"0365{#1}\fi}

 \def\dots{\relax{}
  \ifmmode\def\thedots{\mdots@}\else\def\thedots{\tdots@}\fi %
  \thedots}

 \let\@ldeqno\eqno\let\@ldleqno\leqno
 \def\eqno{\everymath{}\@ldeqno} \def\leqno{\everymath{}\@ldleqno}

  \let\@ldeqalignno\eqalignno
  \def\eqalignno#1{\sam\@ldeqalignno{#1}\mas}
  \let\@ldeqalign\eqalign
  \def\eqalign#1{\sam\@ldeqalign{#1}\mas}

 \def\overrightarrow#1{\vbox{\m@th\ialign{##\crcr
      \rightarrowfill\crcr\noalign{\kern-\p@\nointerlineskip}
      $\hfil\displaystyle{#1}\hfil$\crcr}}}
 \def\overleftarrow#1{\vbox{\m@th\ialign{##\crcr
      \leftarrowfill\crcr\noalign{\kern-\p@\nointerlineskip}
      $\hfil\displaystyle{#1}\hfil$\crcr}}}
 \def\overbrace#1{\mathop{\vbox{\m@th\ialign{##\crcr\noalign{\kern3\p@}
      \downbracefill\crcr\noalign{\kern3\p@\nointerlineskip}
      $\hfil\displaystyle{#1}\hfil$\crcr}}}\limits}
 \def\underbrace#1{\mathop{\vtop{\m@th\ialign{##\crcr
      $\hfil\displaystyle{#1}\hfil$\crcr\noalign{\kern3\p@\nointerlineskip}
      \upbracefill\crcr\noalign{\kern3\p@}}}}\limits}

  \let\@ldmatrix\matrix
  \let\end@ldmatrix\endmatrix
  \def\matrix{\sam\@ldmatrix}
  \def\endmatrix{\end@ldmatrix\mas}
  \let\@ldgather\gather
  \let\end@ldgather\endgather
  \def\gather{\sam\@ldgather}
  \def\endgather{\end@ldgather\mas}
  \let\@ldalign\align
  \let\end@ldalign\endalign
  \def\align{\sam\@ldalign}
  \def\endalign{\end@ldalign\mas}
  \let\@ldaligned\aligned
  \let\end@ldaligned\endaligned
  \def\aligned{\sam\@ldaligned}
  \def\endaligned{\end@ldaligned\mas}
  \let\@ldtag\tag
  \def\tag{\sam\@ldtag}
   %

   \let\MinCDArrowWidth\minCDaw@




\newskip\insertskipamount\newskip\inserthardskipamount
\insertskipamount 6pt plus2pt 
\inserthardskipamount 6pt
\def\insertskip{\vskip\insertskipamount}
\newcount\SplitTest
\def\SetSplitTest{\SplitTest\insertpenalties
  \insert\topins{\floatingpenalty1}%
  \advance\SplitTest-\insertpenalties}
\def\midinsert{\par
 \SaveLastSkip\penalty-150\SetSplitTest\RestoreLastSkip
 \ifnum\SplitTest=-1
  \@midfalse\p@gefalse\else\@midtrue\fi\@ins}
\def\@ins{\par\begingroup\setbox\z@\vbox\bgroup%
  \vglue\inserthardskipamount}
\def\endinsert{\egroup 
  \if@mid \dimen@\ht\z@ \advance\dimen@\dp\z@
    \advance\dimen@\insertskipamount
    \advance\dimen@\pagetotal\advance\dimen@-\pageshrink
    \ifdim\dimen@>\pagegoal\@midfalse\p@gefalse\fi\fi
  \if@mid%
    \ifdim\lastskip<\insertskipamount\removelastskip\insertskip\fi
    \nointerlineskip\box\z@\penalty-200\insertskip
  \else%
    \SaveLastSkip
    \insert\topins{\penalty100 
    \splittopskip\z@skip
    \splitmaxdepth\maxdimen \floatingpenalty\z@
    \ifp@ge \dimen@\dp\z@
    \vbox to\vsize{\unvbox\z@\kern-\dimen@}
    \else \box\z@\nobreak\insertskip\fi}
    \RestoreLastSkip
   \fi\endgroup}


  \newcount\notenumber
  
  \def\note{\advance\notenumber by 1
    \footnote{\the\notenumber)}}

  \newbox\footbox

  \def\footnote#1{\let\@sf\empty
    \ifhmode\edef\@sf{\spacefactor\the\spacefactor}\/\fi
    \sam${}^{\fam0 #1}$\@sf\vfootnote{#1}}%

  \def\vfootnote#1{\insert\footins\bgroup
     \interlinepenalty100 \splittopskip=1pt
     \floatingpenalty=20000
     \leftskip=0pt\rightskip=0pt%
     \parindent=.3em
     \Smallfonts\rm
     \FootItem@{#1}
     \futurelet\next\fo@t}

  \def\FootItem@#1{\par\hangafter1\hangindent=\FootHang
     \setbox0=\hbox{\ignorespaces#1\unskip}%
     \dimen0=.4em\SetOverhang@
     \noindent\rlap{\box0}\kern\Overhang\ignorespaces}


  \def\fo@t{\ifcat\bgroup\noexpand\next \let\next\f@@t
    \else\let\next\f@t\fi \next}
  \def\f@@t{\bgroup\aftergroup\@foot\let\next}
  \def\f@t#1{\baselineskip=10pt\lineskip=1pt
            \lineskiplimit=0pt #1\@foot}%
  \def\@foot{
        \hbox{\vrule height0pt depth5pt width0pt}
        \egroup}
  \skip\footins=12 pt plus 0pt minus 0pt 
  \count\footins=1000 
  \dimen\footins=8in 



 \def\osumess#1{\EdSpider{\immediate\write16{Line \the\inputlineno: #1}}}%
 \def\HideEdStuff{\gdef\EdSpider##1{}}

 \font\BigSym=cmmi10 scaled \magstep 4

 \def\change{\InLMargin{\hbox{\BigSym \char63\kern10pt}}}

 \def\beginchange{\InLMargin{\hbox{\sam\twelvepoint$\heartsuit$\kern10pt}}}

 \def\endchange{\InLMargin{\hbox{\sam\twelvepoint$\spadesuit$\kern10pt}}}

 \def\InLMargin#1{\strut\vadjust{%
     \kern-\strutdepth
     \vtop to \strutdepth{%
         \baselineskip\strutdepth
         \llap{\sam$\smash{\hbox{\EdSpider{#1}}}$}\null}}}

 \def\strutdepth{\dp\strutbox}
 \def\strutheight{\ht\strutbox}

 \def\NoteInRMargin#1{\strut\vadjust{%
     \kern-1.001\strutdepth
     \vtop to \strutdepth{%
       \baselineskip\strutdepth
       \vss\rlap{\ninepoint\unskip\hskip\hsize
         \vtop to 0pt{%
           \hsize=16em\hfuzz=\hsize
           \leftskip=10pt%
           \rightskip=0pt plus 10000pt%
           \baselineskip=9.8pt\lineskip=.2pt%
           \let\\\break
           \noindent\EdSpider{#1}\vss}%
                \kern10pt}\hbox{}}
       }}

 \def\ednote#1{\NoteInRMargin{\tentt #1}}

 \def\cbar{\InLMargin{%
      \dimen0=\strutdepth\advance\dimen0 by \lineskip
      \vrule width 3pt
      height \strutheight depth \dimen0 \kern
      3pt}}

 \def\ccbar{\InLMargin{%
      \dimen0=2\strutdepth\advance\dimen0 by 2\lineskip
      \vrule width 3pt
        height 3\strutheight depth \dimen0 \kern
      3pt}}

 \newinsert\TRMargIns
 \dimen\TRMargIns=\maxdimen

  \def\Ednote#1{\insert\TRMargIns{%
       \vbox to 0pt{\hsize=140pt\hfuzz=\hsize
           \leftskip=6pt%
           \rightskip=0pt plus 10000pt%
           \baselineskip=9.8pt\lineskip=.2pt%
           \let\\\break
           \SetPageRemainder
           \vglue540pt\vglue-\PageRemainder
           \noindent\EdSpider{\tentt #1}\vss}%
       \smallskip}}

 \def\KillEdStuff{\def\ednote##1{}\def\Ednote##1{}%
      \let\change\relax\let\beginchange\relax\let\endchange\relax
       \let\cbar\relax\let\ccbar\relax}


  \topskip=12pt
  \newskip\StdBaselineskip 
  \StdBaselineskip 12pt
  \lineskip=1.1pt
  \lineskiplimit=.8pt
  \widowpenalty=10000 
  \clubpenalty=10000  
  \abovedisplayskip=6pt plus 1pt minus 1pt
  \abovedisplayshortskip=3pt plus 1.5pt
  \belowdisplayskip=6pt plus 1pt minus 1pt
  \belowdisplayshortskip=5pt plus 1pt minus 1pt
  \hfuzz=1.5pt   

  \def\StdPretolerance{100}
  \tolerance=\StdPretolerance

  \newdimen\StdMathsurround
  \StdMathsurround=1.5pt 
  \mathsurround=\StdMathsurround
  \Mas                   

   \def\prose{\relax\hbox{\kern.6\StdMathsurround}}
  
  \def\StdParskip{0pt}    
  \parskip=\StdParskip
  \parindent=0.5cm
 

  \def\Times{ptmr  } 
  \def\TimesI{ptmri  } 
  \def\TimesB{ptmb  }
  \def\TimesBI{ptmbi  }
  \def\HelveticaN{phvrrn }

  =\Times at 10bp
  =\TimesB at 10bp
  \font\tenit=\TimesI at 10bp
  =\TimesBI at 10bp

  \font\tenmrm=cmr10  


    =\Times at 9bp 
    \font\nineit=\TimesI at 9bp 
    =\TimesB at 9bp 
    =\TimesBI at 9bp 

    =\HelveticaN at 9bp 


  =\Times at 12bp
  \font\twelveit=\TimesI at 12bp
  =\TimesB at 12bp


  \font\titleit=\TimesI at 14.4bp
  =\TimesB at 14.4bp

 \SetAuthorHead{AuthorHead} 
 \SetTitleHead{TitleHead}  


  \def\lBr{\raise.125ex\hbox{[\kern.1125ex}}
  \def\rBr{\raise.125ex\hbox{\kern.1125ex]}}

 \setbox\footbox=\hbox{\Smallfonts 2)~}



  \bgroup
  \catcode`\@=11 
  \gdef\itSpacing{%
     \xspaceskip=.31em plus.1em minus.05em \sfcode `f=2001
     \itWarning@\let\itWarning@\itWarning@@}
  \gdef\itSpacingOff{%
     \xspaceskip=0pt \sfcode `f=1000
     \let\itWarning@\relax}
   \global\let\itWarning@\relax
  \gdef\itWarning@@{\errmessage{%
  Special italic spacing already in force
  (you have probably omitted an ``endth'').
  See itSpacing macro in osuPSfnt.sty
         }}
  \egroup

 \fontdimen1\titlebf=0.0pt
 \fontdimen2\titlebf=3.6135pt
 \fontdimen3\titlebf=2.8908pt
 \fontdimen4\titlebf=1.44539pt
 \fontdimen5\titlebf=6.64882pt
 \fontdimen6\titlebf=14.45398pt
 \fontdimen7\titlebf=1.60439pt

 \fontdimen1\tenbi=0.26794pt
 \fontdimen2\tenbi=2.50937pt
 \fontdimen3\tenbi=2.00749pt
 \fontdimen4\tenbi=1.00374pt
 \fontdimen5\tenbi=4.59717pt
 \fontdimen6\tenbi=10.03749pt
 \fontdimen7\tenbi=1.11415pt

 \fontdimen1\twelverm=0.0pt
 \fontdimen2\twelverm=3.01125pt
 \fontdimen3\twelverm=2.409pt
 \fontdimen4\twelverm=1.2045pt
 \fontdimen5\twelverm=5.39615pt
 \fontdimen6\twelverm=12.045pt
 \fontdimen7\twelverm=1.33699pt

 \fontdimen1\twelveit=0.27731pt
 \fontdimen2\twelveit=3.01125pt
 \fontdimen3\twelveit=2.409pt
 \fontdimen4\twelveit=1.2045pt
 \fontdimen5\twelveit=5.37207pt
 \fontdimen6\twelveit=12.045pt
 \fontdimen7\twelveit=1.33699pt

 \fontdimen1\twelvebf=0.0pt
 \fontdimen2\twelvebf=3.01125pt
 \fontdimen3\twelvebf=2.409pt
 \fontdimen4\twelvebf=1.2045pt
 \fontdimen5\twelvebf=5.5407pt
 \fontdimen6\twelvebf=12.045pt
 \fontdimen7\twelvebf=1.33699pt

 \fontdimen1\tenrm=0.0pt
 \fontdimen2\tenrm=2.50937pt
 \fontdimen3\tenrm=2.00749pt
 \fontdimen4\tenrm=1.00374pt
 \fontdimen5\tenrm=4.49678pt
 \fontdimen6\tenrm=10.03749pt
 \fontdimen7\tenrm=1.11415pt

 \fontdimen1\tenit=0.27731pt
 \fontdimen2\tenit=2.50937pt
 \fontdimen3\tenit=2.00749pt
 \fontdimen4\tenit=1.00374pt
 \fontdimen5\tenit=4.47672pt
 \fontdimen6\tenit=10.03749pt
 \fontdimen7\tenit=1.11415pt

 \fontdimen1\tenbf=0.0pt
 \fontdimen2\tenbf=2.50937pt
 \fontdimen3\tenbf=2.00749pt
 \fontdimen4\tenbf=1.00374pt
 \fontdimen5\tenbf=4.61723pt
 \fontdimen6\tenbf=10.03749pt
 \fontdimen7\tenbf=1.11415pt

 \fontdimen1\ninerm=0.0pt
 \fontdimen2\ninerm=2.25842pt
 \fontdimen3\ninerm=1.80673pt
 \fontdimen4\ninerm=0.90337pt
 \fontdimen5\ninerm=4.0471pt
 \fontdimen6\ninerm=9.03374pt
 \fontdimen7\ninerm=1.00273pt

 \fontdimen1\nineit=0.27731pt
 \fontdimen2\nineit=2.25842pt
 \fontdimen3\nineit=1.80673pt
 \fontdimen4\nineit=0.90337pt
 \fontdimen5\nineit=4.02904pt
 \fontdimen6\nineit=9.03374pt
 \fontdimen7\nineit=1.00273pt

 \fontdimen1\ninebf=0.0pt
 \fontdimen2\ninebf=2.25842pt
 \fontdimen3\ninebf=1.80673pt
 \fontdimen4\ninebf=0.90337pt
 \fontdimen5\ninebf=4.15552pt
 \fontdimen6\ninebf=9.03374pt
 \fontdimen7\ninebf=1.00273pt


 \newcount\MaxSpaceFactor
 \MaxSpaceFactor=3000 

 \def\ItemStyle{\rm}
 \def\NrStyle{\rm}
 \def\ItemItemStyle{\rm}

 \MaxItemTag{(iii)}
 \MaxItemItemTag{(iii)}
 \MaxNrTag{(2)}
 \MaxFootTag{2)}
 \def\ReferenceHang{30pt}

 \catcode`\@=\active


\loadbold

=\Times  
=\Times scaled750
=\Times scaled650
\font\rms=\Times scaled 920 

=\TimesBI scaled 860
=\TimesI scaled 860

\textfont0=\rrm  
\scriptfont0=\erm 
\scriptscriptfont0=\srm

\def\Augment#1#2{%
    \toks0\expandafter{#1}\toks2{#2}%
    \edef#1{\the\toks0\the\toks2}}

 \font\twelverma=\Times  scaled 1200
 \font\tenrma=\Times  scaled 1000
 \font\ninerma=\Times scaled 920
 =\Times scaled 840
 \font\sevenrma=\Times scaled 760
 =\Times scaled 680
 \font\fiverma=\Times scaled 600

 \Augment\tenpoint{%
  \textfont0=\tenrma  \scriptfont0=\sevenrma  
  \scriptscriptfont0=\fiverma  }

 \Augment\ninepoint{%
  \textfont0=\ninerma  \scriptfont0=\sevenrma 
  \scriptscriptfont0=\fiverma}

 \Augment\twelvepoint{%
  \textfont0=\twelverma  \scriptfont0=\ninerma  
  \scriptscriptfont0=\sevenrma}

\mathsurround=1pt
\hsize=13.45truecm
\vsize=19.5truecm
\hoffset=1.25truecm
\voffset=2truecm
\advance\baselineskip by 2pt

\predefine\til{\~}
\def\~#1{\relax\ifmmode\widetilde{#1}\else\til{#1}\fi}

\redefine \le{\leqslant}
\redefine \ge{\geqslant}
\define \wt#1{\mathaccent"0365{#1}}
\define \wh#1{\mathaccent"0362{#1}}

\define \coo{\Cal O_0}

\define \iss{\,\Mathaccent{\raise -.8 ex\hbox{$\widetilde{}$\kern.1em}}\rightarrow\,}

\define \expp{\operatorname{\fam0 exp}}

\define \tpp{\mathop{\fam0 top}}
\define \ab{\mathop{\fam0 ab}}

\define \lln{\operatorname{\fam0 log\,}}

\define \kr{\mathop{\fam0 ker}}
\define \Frob{\operatorname{\fam0 Frob}}
\define \chr{\mathop{\fam0 char}\,}

\define \res{\operatorname{\fam0 res}}

\define \Tr{\operatorname{\fam0 Tr\,}}

\define \Gal{\mathop{\fam0 Gal}}

\Mas
\HideEdStuff
\rm 
 

\def\issn{{\nineit ISSN 1464-8997 (on line) 1464-8989 (printed)}}

\def\gtp{{\nineit Published 10 December 2000: \ \copyright\ Geometry \& 
Topology Publications}}

\def\gtv3{{\nineit Geometry \& Topology Monographs, Volume 3 (2000) --
Invitation to higher local fields}}


\def\lione
{{\rms Geometry \& Topology Monographs}}

\def \litwo{{\rms Volume 3: Invitation to higher local fields
}} 

\def\tinfo #1.#2.#3-#4
{{
\noindent  {\lione} \hfill 
\par 
\vskip-1.5pt
\noindent {\litwo} \hfill
\par 
\vskip-1,5pt
\noindent {\rms Part #1, section #2, pages #3--#4} \hfill
\vskip24pt 
}}

\def\tinfos #1.#2.#3-#4
{{
\noindent  {\lione} \hfill 
\par 
\vskip-1.5pt
\noindent {\litwo} \hfill
\par 
\vskip-1.5pt
\noindent {\rms Pages #3--#4} \hfill
\vskip24pt 
}}

\def\tinfoi #1
{{
\noindent  {\lione} \hfill 
\par 
\vskip-1.5pt
\noindent {\litwo} \hfill
\par 
\vskip-1.5pt
\noindent {\rms Pages iii--xi: Introduction and contents} \hfill
\vskip26pt 
}}


  \def\titlepagehead{\hfil}

  \newif\iftitlepage\titlepagefalse
  \newif\ifblankpage\blankpagefalse
  \def\makeheadline{
     \ifblankpage{}\else%
     \iftitlepage
\vbox{\line{\vbox to 8.5pt{}
\ninerm
\copy\HLinebox \hfill
\hglue5mm\ninebf\folio 
\titlepagehead}}%
      \else
\vbox{\ifodd\pageno\rightheadline\else\leftheadline\fi}%
      \fi\vskip 12pt\fi}%
     \def\rightheadline{\line{\vbox to 8.5pt{}%
      \ninerm
\copy\TitleBox \hfill
\hglue5mm\ninebf\folio}}%
     \def\leftheadline{\line{\vbox to 8.5pt{}%
        \unskip\ninerm\unskip\ninebf\folio\hglue5mm
 \hfill \copy\AuthorBox
}}

 \footline={\ifblankpage{}\else
\iftitlepage\ninepoint\sam\hfill
\line{\vbox to 8.5pt{}
\copy\TFLinebox
\hfill
\hglue5mm 
}
            \else
\ninepoint\sam\hfill
\line{\vbox to 8.5pt{}
\copy\FLinebox
\hfill 
\hglue5mm
}
\hfil\fi\global\titlepagefalse\fi}

\def\blankpage{{\blankpagetrue\noindent\hbox to 10pt{\hss}\vfill
\pagebreak}}

\tenpoint\rm 
 